\documentclass[12pt]{amsart}

\usepackage{amsmath}
\usepackage{amscd}
\usepackage{amssymb}
\usepackage{latexsym}
\usepackage{stmaryrd}
\usepackage{mathrsfs}
\DeclareMathAlphabet{\mathbbb}{U}{bbold}{m}{n}

\newtheorem{theorem}{Theorem}[section]
\newtheorem*{theorem*}{Theorem}
\newtheorem{lemma}[theorem]{Lemma}
\newtheorem*{lemma*}{Lemma}
\newtheorem{corollary}[theorem]{Corollary}
\newtheorem{proposition}[theorem]{Proposition}

\newtheorem{remark}[theorem]{Remark}
\newtheorem{definition}[theorem]{Definition}

\makeatletter
\def\revddots{\mathinner{\mkern1mu\raise\p@
\vbox{\kern7\p@\hbox{.}}\mkern2mu
\raise4\p@\hbox{.}\mkern2mu\raise7\p@\hbox{.}\mkern1mu}}
\makeatother 
\newcommand{\bgl}{\begin{equation}} 
\newcommand{\egl}{\end{equation}}
\newcommand{\bgloz}{\begin{equation*}} 
\newcommand{\egloz}{\end{equation*}}
\newcommand{\bgln}{\begin{eqnarray}} 
\newcommand{\egln}{\end{eqnarray}}
\newcommand{\bglnoz}{\begin{eqnarray*}} 
\newcommand{\eglnoz}{\end{eqnarray*}}
\newcommand{\btheo}{\begin{theorem}}
\newcommand{\etheo}{\end{theorem}}
\newcommand{\btheooz}{\begin{theorem*}}
\newcommand{\etheooz}{\end{theorem*}}
\newcommand{\blemma}{\begin{lemma}}
\newcommand{\elemma}{\end{lemma}}
\newcommand{\blemmaoz}{\begin{lemma*}}
\newcommand{\elemmaoz}{\end{lemma*}}
\newcommand{\bproof}{\begin{proof}}
\newcommand{\eproof}{\end{proof}}
\newcommand{\bbew}{\begin{beweis}}
\newcommand{\ebew}{\end{beweis}}
\newcommand{\bremark}{\begin{remark}\em}
\newcommand{\eremark}{\end{remark}}
\newcommand{\bdefin}{\begin{definition}}
\newcommand{\edefin}{\end{definition}}
\newcommand{\bprop}{\begin{proposition}}
\newcommand{\eprop}{\end{proposition}}
\newcommand{\bcor}{\begin{corollary}}
\newcommand{\ecor}{\end{corollary}}
\newcommand{\bfa}{\begin{cases}} 
\newcommand{\efa}{\end{cases}}
\newcommand{\cD}{\mathcal D}

\newcommand{\cQ}{\mathcal Q}

\def\Az{\mathbb{A}}

\def\Nz{\mathbb{N}}

\def\Qz{\mathbb{Q}}

\def\Zz{\mathbb{Z}}

\def\1z{\mathbbb{1}}
\newcommand{\fA}{\mathfrak A}

\newcommand{\an}[1]{``#1''} 

\newcommand{\ri}{\rightarrow}

\newcommand{\ma}{\mapsto} 
\newcommand\into{\hookrightarrow} 

\def\SEMI{\mbox{$\times\kern-2pt\vrule height5pt width.6pt \kern3pt $}}

\newcommand{\End}{{\rm End}\,}

\newcommand{\Ad}{{\rm Ad\,}}

\newcommand{\reg}{^\times} 
\newcommand{\lspan}{{\rm span}} 
\newcommand{\clspan}{\overline{\lspan}} 

\newcommand{\fa}{\text{ for all }} 
\newcommand{\ilim}{\varinjlim} 
\newcommand{\rte}{\overset{e}{\rtimes}} 
\newcommand{\lge}{\left\{} 
\newcommand{\rge}{\right\}} 
\newcommand{\gekl}[1]{\lge #1 \rge} 
\begin{document}

\title[Dilations of semigroup crossed products]{Dilations of semigroup crossed products as crossed
products of dilations}
\author[N. S. Larsen]{Nadia S. Larsen}
\address{Department of Mathematics, University of Oslo, PO BOX 1053 Blindern, N-0316 Oslo, Norway.}
\email{nadiasl@math.uio.no}
\author[X. Li]{Xin Li}
\address{Department of Mathematics, Westf{\"a}lische Wilhelms-Universit{\"a}t M{\"u}nster, Einsteinstra{\ss}e 62, 48149 M{\"u}nster, Germany}
\email{xinli.math@uni-muenster.de}

\thanks{This research was supported by the Research Council of Norway and the Deutsche Forschungsgemeinschaft. The first named author thanks
J. Cuntz and S. Echterhoff for their kind hospitality during a sabbatical visit at Westf{\"a}lische Wilhelms-Universit{\"a}t
M{\"u}nster in October 2009, where this research was initiated. The second named author thanks the operator algebra group in Oslo for a nice visit at the University of Oslo.}

\begin{abstract}
Laca constructed a minimal automorphic dilation for every semigroup dynamical system arising from an action of an Ore semigroup by injective endomorphisms of a unital $C^*$-algebra. Here
we show that the semigroup crossed product with its action by inner endomorphisms given by the implementing isometries
has as minimal automorphic dilation the group crossed product of the original dilation. Applications include
recent examples studied by Cuntz and the second named author.
\end{abstract}
\date{September 21, 2010.}

\maketitle

\section*{Introduction}

Starting from the work of Cuntz in \cite{Cun0}, the theory of semigroup crossed products has
consolidated itself as a rich source of many new and interesting examples
of $C^*$-algebras. In \cite{Cun0}, semigroup crossed products appear as corners in ordinary
group crossed products. Later, in \cite{S}, a definition of semigroup crossed products was introduced in terms
of a universal property for a covariance relation relating representations of the $C^*$-algebra and
representations of the semigroup. Much of the subsequent theory developed along this path,
see for example \cite{M} and \cite{LR1}. However, the corner representation in an
ordinary group crossed product played an important role, see \cite{M}. Laca
extended the work of Murphy to a large class of not necessarily abelian semigroups, the
Ore (right-reversible) semigroups which
act by injective endomorphisms of a unital $C^*$-algebra, see \cite{La}.

Recent new examples of $C^*$-algebras arising in number theory and in purely algebraic
context turned out to be instances of semigroup crossed products, see \cite{Cun1}, \cite{Cun-Li1}, \cite{Cun-Li2} and
\cite{Li}. From these semigroup crossed products isomorphisms emerge involving ordinary group crossed
products by affine-type groups. These isomorphisms play an important role in the $K$-theory computations.

The goal of this short note is to establish a general
result about the minimal automorphic dilation of a certain type of semigroup dynamical system.
Given a semigroup dynamical system, one can first dilate the dynamical system to the effect that the endomorphisms become automorphisms and then form the crossed product (by a group), or one can first form the semigroup crossed product and then dilate the resulting C*-algebra. Our main result (Theorem~\ref{main}) says that the outcome is the same.

We were led to this result by certain isomorphisms established in \cite{Cun1} and \cite{Cun-Li1} (see Section~\ref{App}). They turn out to be special cases of our main theorem. Therefore, this short note gives a unified treatment of these phenomena.

\section{Dilations of semigroup crossed products}
Let $A$ be a unital C*-algebra and
$\alpha$ an action of a semigroup $P$ on $A$ by endomorphisms, so
$\alpha$ is a semigroup homomorphism $P \ri \End(A)$. We associate the semigroup crossed product $A \rte_{\alpha} P$ with the semigroup C*-dynamical system $(A,P,\alpha)$. The symbol $\rte$ stands for \an{semigroup crossed product by endomorphisms}. The $C^*$-algebra $A \rte_{\alpha} P$ is
universal for covariant representations of $(A, P,\alpha)$: these are pairs
$(\pi, W)$ consisting of a non-degenerate representation $\pi$ of $A$ on a Hilbert space $H$ and a homomorphism $W$ from $P$ into the
semigroup of isometries on $H$ such that the
\emph{covariance condition} $\pi(\alpha_p(a))=W_p\pi(a)W_p^*$ is satisfied for
all $a\in A$ and $p\in P$. The crossed product $A \rte_{\alpha} P$ is generated
as a $C^*$-algebra by a universal covariant representation $(i_A, i_P)$, see
for example \cite{S, M, LR1}. One can argue that the idea behind the definition of $A \rte_{\alpha} P$ is to enlarge $A$ by adding isometries
$w_p=i_P(p)$ for all $p \in P$ so that the action of
$\alpha$ becomes inner, in the sense that $\alpha_p(\cdot) = w_p \cdot w_p^* \fa p \in P$.
In general, it may happen that $A \rte_{\alpha} P$ is zero, see e.g. \cite[Example 2.1(a)]{S}.

In this note we will be concerned with actions of Ore semigroups.
Recall that a (right-reversible) Ore semigroup is a cancellative semigroup $P$ such that $Ps\cap Pt\neq \emptyset$ for
all $s,t\in P$. There is a similar notion of left-reversible Ore semigroups.
A cancellative semigroup is Ore precisely when it can be embedded in a group $G$ in such a way
that $G=P^{-1}P$. $G$ is called the enveloping group of $P$ and is determined uniquely up to isomorphism by $P$, cf. \cite[Theorems 1.23, 1.24, 1.25]{Cli-Preston}. The
semigroup structure of $P$ induces a partial pre-order on $P$ (and also on $G$) given by $p\leq r$ if $r\in Pp$.

Let $\alpha$ be an action of an Ore semigroup $P$ on a unital $C^*$-algebra $A$.
If there is a non-trivial covariant representation for $(A, P, \alpha)$ then
$A \rte_{\alpha} P$ is non-zero by e.g. \cite{LR1}, and due to the
right-reversibility of $P$ the crossed product $A \rte_{\alpha} P$ is the closed
span of the monomials $i_P(p)^*i_A(a)i_P(r)$ for $a\in A$ and $p,r\in P$, see \cite{La}.

Fix a semigroup dynamical system $(A, P, \alpha)$ with $P$ an Ore semigroup
and $\alpha$ an action of $P$ by endomorphisms $\alpha_p$ which are injective
for all $p\in P$.
Laca used dilation theory for projective representations of semigroups to show that the existence of a non-trivial
covariant representation for this kind of $(A, P, \alpha)$ is automatic.
Indeed, if $(A, P, \alpha)$
is a semigroup dynamical system with $P$ an Ore semigroup and $\alpha$ an
action by injective endomorphisms, let $G=P^{-1}P$ be the enveloping group of $P$, and recall that
\cite[Theorem 2.1]{La} shows that there exist an ordinary C*-dynamical system
$(B,G,\beta)$ and an embedding
$\iota: A \into B$ such that
\bgl
\label{beta_extends_alpha}
  \beta_p \circ \iota = \iota \circ \alpha_p \fa p \in P
\egl
and
\bgl
\label{dense}
  \overline{\bigcup_{p \in P} \beta_p^{-1}(\iota(A))} = B.
\egl
The triple $(B,G, \beta)$ is the \emph{minimal automorphic dilation} of $(A, P,
\alpha)$ and is unique up to isomorphism. The identities in
\eqref{beta_extends_alpha} express the dilation phenomena, and
condition \eqref{dense} is the minimality. Laca's other main result shows that $A\rte_\alpha P$
can be realised as a full corner of the ordinary crossed product $B\rtimes_\beta G$, cf. \cite[Theorem 2.4]{La}.

\bremark\label{Ainfty_model}
A model for $B$ is given as follows (see the proof of \cite[Theorem 2.1]{La}): Take the direct limit over the directed set $P$ of copies $A_p$ of $A$
for all $p\in P$ with connecting homomorphisms $\alpha_r^p:A_p\to A_r$ given by
$\alpha_r^p:=\alpha_{rp^{-1}}$
when $p\leq r$. Then there is an action
$\alpha_\infty$ of $G$ on
$A_{\infty, \alpha} = \ilim_{P} \gekl{A; \alpha_r^p}_{r\leq p}$ such that $(A_{\infty,
\alpha}, G, \alpha_\infty)$ and the embedding $\iota^1:A_e\to A_{\infty, \alpha}$ satisfy
\eqref{beta_extends_alpha}
and \eqref{dense}.
Here $e$ is the identity element of $P$. We always assume that our semigroups have identity elements.
\eremark

Let $1_A$ denote the unit of $A$, let $w_p=i_P(p)$ for $p\in P$
and let $u_g$, $g\in G$, be the unitaries in the multiplier
algebra of $B \rtimes_{\beta} G$ which implement $\beta$.
The inclusion $\iota: A \into B$ induces an inclusion $A \rte_{\alpha} P \into B
\rtimes_{\beta} G$, also denoted $\iota$, with the property that
\bgl
\label{w_p}
  \iota(w_p) = u_p \iota(1_A)
\egl
for $p\in P$, see \cite[Theorem 2.4]{La}. We note the following

\blemma\label{approx_unit}
The net $(u_p^* \iota(1_A) u_p)_{p \in P}$ is an approximate unit in $B$.
\elemma
\bproof
This is an immediate consequence of the fact that $u_p^*\cdot u_p =
\beta_p^{-1}(\cdot)$ for all $p$ in $P$ and that $\bigcup_{p \in P} u_p^* \iota(A) u_p$ is dense in $B$ by \eqref{dense}.
\eproof

Next we notice that we can form a new semigroup dynamical system with injective
endomorphisms by
letting $\phi$ be the action of $P$ on $A \rte_{\alpha} P$ given by
$\phi_p:=\Ad w_p$ for $p\in P$. Moreover, let $\psi$ be the action of $G$ on $B \rtimes_{\beta} G$ given by $\psi_g := \Ad u_g$ for $g$ in $G$.

The main thrust of this note is that the dilation of the new system $(A \rte_{\alpha} P, P, \phi)$ is
the crossed product $B \rtimes_{\beta} G$. This result
was suggested by very recent examples of actions arising in number theoretic
and purely algebraic context, see \cite{Cun1, Cun-Li1, Cun-Li2}.

\btheo\label{main} The embedding $\iota:A \rte_{\alpha} P \into B \rtimes_{\beta} G$ and the
 dynamical system $(B \rtimes_{\beta} G, G, \psi)$
 form the
minimal automorphic dilation of $(A \rte_{\alpha} P, P, \phi).$
\etheo

\bproof
By \cite[Theorem 2.1]{La} we have to show the dilation phenomena and the
minimality. That $\psi_p$ dilates $\phi_p$ for every $p$ in $P$ in the sense of
\eqref{beta_extends_alpha} follows from \eqref{w_p}.

To show that $(B \rtimes_{\beta} G, G, \psi)$ is minimal, we need to
show that any $y\in B \rtimes_{\beta} G$ can be approximated from
$\bigcup_{p\in P} u_p^* (\iota(A \rte_{\alpha} P)) u_p$. It suffices to
prove this claim for $y=bu_g$ for $b$ in $B$ and $g$ in $G$ because
elements of this form span a dense
subset of $B \rtimes_{\beta} G$. So fix $b$ and $g$.
By \eqref{dense}, we can find a net $(a_p)_{p \in P}$ in $A$ such that
\bgl
\label{b}
  b = \lim_{p} u_p^* \iota(a_p) u_p.
\egl
Since $G = P^{-1} P$, we can find for all $p$ in $P$ elements $q_p$ and $r_p$ in $P$
such that $q_p^{-1} r_p = p g p^{-1}$. Then $u_{p g p^{-1}}=u_{q_p^{-1} r_p}=
u_{q_p}^* u_{r_p}$. Hence, using \eqref{w_p}, we obtain
\begin{align}
u_p^* (\iota(a_p w_{q_p}^* w_{r_p})) u_p
&=u_p^* \iota(a_p) \iota(w_{q_p}^*) \iota(w_{r_p}) u_p\notag\\
&=u_p^* \iota(a_p)u_{q_p}^* u_{r_p} \iota(1_A) u_p\notag\\
&=u_p^* \iota(a_p) u_p u_g u_p^* \iota(1_A) u_p.\label{image_psi_p}
\end{align}
Lemma~\ref{approx_unit} and \eqref{b} imply that the term in \eqref{image_psi_p}
converges to $bu_g$. Hence $b u_g$ can be approximated from $\bigcup_{p\in P}
u_p^* (\iota(A \rte_{\alpha} P)) u_p$, as claimed.
\eproof

\bremark
In other words, we have shown that the C*-algebras
\bgloz
  B \rtimes_{\beta} G \text{ and } \ilim_{p \in P} \gekl{A \rte_{\alpha} P; \phi_p}
\egloz
are canonically isomorphic.
\eremark

In the cases of interest for our applications the semigroup $P$ has the form $P=
H\rtimes_\eta M$ where $M$ is a right-reversible Ore semigroup, $H$ is a group and
$\eta:M\to \End(H)$ an action of $M$ by injective endomorphisms
of $H$. Both $H$ and $M$ can be embedded as subsemigroups in $P$ via $H \ni h \ma (h, e_M) \in P$ and $M \ni m \ma (e_H, m) \in P$, where $e_H$ and $e_M$ denote the identity elements in $H$ and $M$, respectively. Since $H$ is a group, the action $\alpha$ of $P$ restricts to an action $\alpha \vert_H$ of $H$ by automorphisms (as long as the identity element of $P$ acts as the identity homomorphism). It follows from our assumptions that $P$ is again an Ore semigroup. Let $G$ be the enveloping group of $P$. The equality $(-h,e_M)(h,m)=(e_H,m)$
shows that the semigroup $M$ is a cofinal subset of $P$ (i.e. for any
$(h,m)\in P$ there is an element $x\in M$ such that $(h,m)\leq x$ in the right-order on $P$). Therefore
an inductive limit over $P$ can be viewed as an inductive limit over $M$.
Hence Theorem~\ref{main} implies that the
embedding $\iota:A\rte_{\alpha} P \into B\rtimes_\beta G$
and the triple $(B\rtimes_\beta G, M^{-1}M, \psi \vert_{M^{-1}M})$ form the minimal automorphic dilation of the semigroup dynamical system
$(A\rte_{\alpha} P, M, \phi \vert_M)$. With the notation of
Remark~\ref{Ainfty_model},
\bgl\label{dilation_absorbs_inner}
(A\rte_{\alpha} P)_{\infty, \phi \vert_M}\cong B\rtimes_\beta G.
\egl

\bremark Note that in the case $P=H\rtimes_\eta M$, the condition \eqref{dense} expressing the minimality of the minimal automorphic dilation $(B,G,\beta)$ takes the following form:
\bgloz
  \bigcup_{m \in M}\beta_{(e_H, m)}^{-1} (\iota(A)) \text{ is dense in } B.
\egloz
This follows from cofinality of $M$ in $P$.
\eremark

The motivating example is that of
$\Zz\rtimes \Nz^\times$, where the multiplicative semigroup $\Nz^\times$ of
non-zero natural numbers acts on $\Zz$ by multiplication $n \mapsto kn$ for $k \in \Nz^\times$ and
$n \in \Zz$. It is easily seen that $\Zz\rtimes\Nz^\times$ is right-reversible. Note that already $\Nz\rtimes \Nz^\times$
is right-reversible, where $\Nz$ is the additive semigroup of natural numbers with $0$. But $\Nz^\times$ is not cofinal in
$\Nz\rtimes\Nz^\times$ for the right-order. However, $\Nz\rtimes\Nz^\times$ is not left-reversible, since
for pairs $(l,k),(n,m)\in \Nz\rtimes\Nz^\times$ such that $n-l$ is not congruent to $0$ mod $\operatorname{gcd}(k,m)$ we have
$(l,k)(\Nz\rtimes\Nz^\times) \cap (n,m)(\Nz\rtimes\Nz^\times) = \emptyset$, see \cite[Proposition 2.2]{La-Ra2}
(there it is shown that $(\Qz \rtimes \Qz \reg_+, \Nz\rtimes\Nz^\times)$ is quasi-lattice ordered with respect to its left-order).

\section{Applications}
\label{App}

\subsection{The $C^*$-algebra $\cQ_\Nz$}\label{first-subsection} The first application concerns the
$C^*$-algebra $\cQ_\Nz$ constructed in \cite{Cun1}. Recall from \cite{Cun1}
that $\cQ_\Nz$ is the universal $C^*$-algebra generated by isometries
$s_m$ for $m\in \Nz^\times$ and a unitary $u$ subject to the relations $s_k s_m = s_{km}$,
$s_mu=u^ms_m$ and $\sum_{n=0}^{m-1}u^ns_ms_m^*u^{-n}=1$ for $k, m\in \Nz^\times$. Let $P_\Nz$
be the semidirect product $\Nz\rtimes \Nz^\times$ arising from the action of $\Nz^\times$
on $\Nz$ given by $k \cdot n=kn$ for $k\in \Nz^\times$ and $n\in \Nz$. (Note that this $P_\Nz$ is
the opposite semigroup of the semigroup of matrices denoted by the same symbol in \cite{Cun1}.)
The operation in $P_\Nz$ is given by $(l, k)(n, m)=(l+kn, km)$ for $(l, k), (n, m)\in P_\Nz$.

Let $e_m=s_ms_m^*$ denote the range projection of $s_m$ for all $m \in \Nz^\times$.
It is shown in \cite{Cun1} that the $C^*$-algebra $\cD$ generated
by the projections $u^ne_mu^{-n}$ where $0\leq n<m$ for all $m\in \Nz^\times$ is commutative. The semigroup
$P_\Nz$ acts on $\cD$ by $\alpha_{(n,m)}:= \Ad (u^ns_m)$. Let $(i_{\cD}, i_{P})$ denote the
universal covariant pair of the semigroup system $(\cD, P_\Nz, \alpha)$. Let $w_{(n, m)}
=i_{P}(n,m)$ for $n\in \Zz$, $m\in \Nz^\times$. It follows that the isometries
$w_{(0,m)}$ for $m\in \Nz^\times$ and the unitary $w_{(1,1)}$ satisfy the relations defining
$\cQ_\Nz$. By the universal property of $\cQ_\Nz$ there is a homomorphism from  $\cQ_\Nz$
to $\cD\rte_\alpha P_\Nz$ sending $u$ to $w_{(1,1)}$ and $s_m$ to $w_{(0,m)}$ for all $m\in \Nz^\times$.
The algebra $\cD\rte_\alpha P_\Nz$ is the closed span of monomials of the form
$w_{(l, k)}^* i_{\cD}(f)w_{(n, m)}$ for $f\in \cD$; $(l, k), (n, m)\in P_\Nz$. Writing
$w_{(n,m)}=w_{(1,1)}^nw_{(0,m)}$ and using that $\cD$ is generated by projections
of the form $u^ne_mu^{-n}$  shows that $\cD\rte_\alpha P_\Nz$ is generated as a $C^*$-algebra
by the elements $w_{(1,1)}$ and $w_{(0,m)}$ for $m\in \Nz^\times$. Simplicity of
$\cQ_\Nz$ implies that the homomorphism defined above is an isomorphism, i.e.
\begin{equation}\label{Qn_as_crossed_product_by_axplusb}
  \cQ_\Nz\cong \cD \rte_\alpha P_\Nz.
\end{equation}
This fact is implicit in \cite{Cun1}, and follows from taking $R=\Zz$ in \cite{Cun-Li1}, see Paragraph \ref{second-subsection}.

Since $u$ is a unitary, $\alpha$ extends to an action of
$\Zz\rtimes \Nz^\times$ on $\cD$ (also denoted by $\alpha$).
Moreover, $\cD \rte_\alpha P_\Nz$ and $\cD \rte_\alpha (\Zz \rtimes \Nz \reg)$ are canonically isomorphic.
It is easy to see that the $ax+b$-group $P_\Qz^+=\Qz \rtimes \Qz \reg _+$ is the enveloping group
of $\Zz \rtimes \Nz \reg$.
By Theorem~\ref{main} as formulated in \eqref{dilation_absorbs_inner},
$(\cD\rte_\alpha P_\Nz)_{\infty,\phi \vert_{\Nz^\times}}$
is isomorphic to $\cD_{\infty, \alpha} \rtimes_\beta P_\Qz^+$.

It is proved in \cite{Cun1} that $\cD$ is isomorphic to $C(\hat{\Zz})$; the isomorphism
sends $u^ne_mu^{-n}$ to the characteristic function of the set $n+m\hat{\Zz}$, see \cite{Cun-Li1}.
Here $\hat{\Zz}$ is the profinite completion of the integers. Under this
isomorphism the action $\alpha$ of $P_\Nz$ on $\cD$ is transformed into
the action of $P_\Nz$ on $C(\hat{\Zz})$ given by affine transformations,
\bgl
  \alpha_{(n,m)}^{\operatorname{aff}}(f)(x) =
    \begin{cases}f(m^{-1}(x-n)),&\text{ if } x \in n+m\hat{\Zz}\\
    0, \text{ otherwise}.
    \end{cases}
\egl

It follows from \cite{Cun1} and \cite{Cun-Li1} that $\cD_{\infty, \alpha}$ is isomorphic to $C_0(\Az_f)$,
where $\Az_f$ is the finite adele ring over the rationals. The
embedding of $\cD$ into $\cD_{\infty, \alpha}$ is the canonical embedding
$\iota:C(\hat{\Zz})\to C_0(\Az_f)$ obtained from cutting down with the projection equal to the characteristic function of
$\hat{\Zz}$. The action $\beta$ of $P_\Qz^+$ extending $\alpha^{\operatorname{aff}}$ is the
action by affine transformations on $C_0(\Az_f)$. Taking into account
\eqref{Qn_as_crossed_product_by_axplusb}, we have established that
$C_0(\Az_f)\rtimes_{\alpha^{\operatorname{aff}}} P_\Qz^+$ is isomorphic to
the dilation of $\cQ_\Nz$ by the endomorphisms $x \mapsto s_m x s_m^*$, $m\in \Nz^\times$.
This dilation is the algebra $\bar{\cQ}_\Nz$ from \cite{Cun1}, and we have recovered \cite[Theorem 6.4]{Cun1}.

\subsection{$C^*$-algebras associated to integral domains}
\label{second-subsection}
Similarly to the previous example we now outline how to use
Theorem~\ref{main} to recover \cite[Theorem 2]{Cun-Li1}.

Recall the following from \cite{Cun-Li1}: Suppose that
$R$ is an integral domain with set of units strictly contained in
$R^\times=R\setminus \{0\}$ and such that for each
$m\in R^\times$ the ideal $(m)$ of $R$ is of finite index in $R$. Then
$\fA[R]$ is the universal $C^*$-algebra generated by isometries
$\{s_m\mid m \in R^\times\}$ and unitaries $\{u^n \mid n \in R\}$
subject to a number of relations (in the case $R=\Zz$ these are the relations defining $\cQ_\Nz$).
The $C^*$-algebra $\cD[R]:=\clspan\{u^ne_mu^{-n}\mid n\in R, m\in R^\times\}$ is
commutative, and is isomorphic to the algebra of continous functions
on the profinite completion $\hat{R}$ of $R$.

To apply Theorem~\ref{main} we need to start with a semigroup crossed product, and
dilate it to obtain an ordinary group crossed product. Therefore we proceed
in the opposite direction as that pursued in \cite{Cun-Li1} and start by
identifying $\fA[R]$ with a semigroup crossed product. Indeed, let $P_R$ be the
right-reversible semigroup $R \rtimes R \reg$ where $R \reg$ acts multiplicatively on the additive group $R$.
Let $\alpha_{(l,k)}:=\Ad (u^ls_k)$ for $(l,k)\in P_R$. Then $\alpha_{(l,k)}(u^ne_mu^{-n})=u^{l+kn}e_{km}u^{-(l+kn)}$
for all $(l, k),(n,m)\in P_R$, and it follows that $\alpha$ is an action
of $P_R$ by injective endomorphisms of $\cD[R]$. If $(i_\cD, i_P)$
is the universal covariant pair of $(\cD[R], P_R, \alpha)$, then $i_\cD(1)$
is the identity element in the semigroup crossed product and each
$w_{(n,m)}:=i_P(n,m)$ is an isometry for $(n,m)\in P_R$.
It follows from our construction that the isometries $w_{(0,m)}$ for $m\in R^\times$ and the
unitaries $w_{(n,1)}$ for $n \in R$ satisfy the relations defining $\fA[R]$.
Moreover, it follows as in Section~\ref{first-subsection} that
$w_{(n,1)}$ and $w_{(0,m)}$ for $n \in R$ and $m \in R \reg$ generate $\cD[R] \rte_\alpha P_R$.
Then simplicity of $\fA[R]$, established in \cite{Cun-Li1}, implies that there is an isomorphism of
$\fA[R]$ onto $\cD[R]\rte_\alpha P_R$ which carries $s_m$ to $w_{(0,m)}$ for $m\in R^\times$.

The enveloping group of $P_R$ is $P_{Q(R)}$, the $ax+b$-group over the quotient field $Q(R)$ of $R$.
Theorem~\ref{main} implies that $(\cD[R])_{\infty, \alpha} \rtimes_\beta P_{Q(R)}$ is isomorphic to the dilation
$(\cD[R]\rte_\alpha P_R)_{\infty, \phi \vert_{R^\times}}$, and so there is an isomorphism
\bgl\label{dilation_of_AR} (\cD[R])_{\infty, \alpha} \rtimes_\beta P_{Q(R)}\cong (\fA[R])_{\infty, \Ad\vert_{R^\times}}.
\egl

It was shown in \cite{Cun-Li1} that $(\cD[R])_{\infty, \alpha}$ is isomorphic to $C_0(\mathscr{R})$
where $\mathscr{R}$ is the generalised finite adele ring associated with $R$.
The $C^*$-algebra $(\fA[R])_{\infty, \Ad\vert_{R^\times}}$ is the stabilisation $\fA(R)$ from \cite{Cun-Li1},
and \eqref{dilation_of_AR} recovers therefore the isomorphism of $\fA(R)$ with $C_0(\mathscr{R})\rtimes P_{Q(R)}$.
This is \cite[Theorem 2]{Cun-Li1}.


\begin{thebibliography}{99}


\bibitem[Cli-Pre]{Cli-Preston} A.H. Clifford and G. B. Preston, The algebraic
theory of semigroups - Vol. I, Mathematical Surveys 7, American Mathematical
Society, Providence, Rhode Island, 1961.

\bibitem[Cu1]{Cun0} J. Cuntz, \emph{Simple $C^*$-algebras generated by isometries}, Comm. Math. Phys. {\bf 57} (1977), 173--185.


\bibitem[Cu2]{Cun1} J. Cuntz,
\emph{$C^*$-algebras associated with the $ax+b$-semigroup over $\Nz$},
in $K$-Theory and noncommutative geometry (Valladolid, 2006),
European Math. Soc., 2008, pp 201--215.


\bibitem[Cu-Li1]{Cun-Li1} J. Cuntz and X. Li, \emph{The regular $C^*$-algebra
of an integral domain},
to appear in the proceedings of the conference in honour of A. Connes' $60$th
birthday.

\bibitem[Cu-Li2]{Cun-Li2} J. Cuntz and X. Li, \emph{$C^*$-algebras associated
with integral domains and crossed products by actions on adele spaces},
J. Noncomm. Geom., to appear.

\bibitem[La]{La} M. {Laca},
  \emph{From endomorphisms to automorphisms and back: dilations and full corners},
  J. London Math. Soc. \emph{61} (2000), 893--904.

\bibitem[La-Ra1]{LR1} M. Laca and I. Raeburn, Semigroup crossed products and the
Toeplitz algebras of nonabelian groups, {\it J. Funct. Anal.} {\bf 139} (1996),
415--440.

\bibitem[La-Ra2]{La-Ra2} M. Laca and I. Raeburn, \emph{Phase transition on the
Toeplitz algebra of the affine semigroup over the natural numbers}, preprint
arXiv:math.OA/0907.3760v2.

\bibitem[Li]{Li} X. Li, \emph{Ring $C^*$-algebras}, Math. Ann., to appear.

\bibitem[Mur]{M} G. J. Murphy, Crossed products of $C^*$-algebras by
endomorphisms, {\it Integral Equations Oper. Theory} {\bf 24} (1996), 298--319.


\bibitem[Sta]{S} P. J. Stacey, Crossed products of $C^*$-algebras by
$*$-endomorphisms, {\it J. Austral. Math. Soc.} Ser. A {\bf 54} (1993), 204--212.

\end{thebibliography}
\end{document}